\documentclass[12pt,reqno]{amsart}
\usepackage{ amsmath , amssymb }
\usepackage[mathscr]{eucal}

\newtheorem{thm}{Theorem}
\newtheorem*{thma}{Theorem A}
\newtheorem*{thmb}{Theorem B}
\newtheorem{conj}{Conjecture}
\newtheorem{cor}[thm]{Corollary}

\newtheorem{lemma}[thm]{Lemma}

\theoremstyle{definition}\newtheorem{rmk}{Remark}
\theoremstyle{definition}
\theoremstyle{definition}
\numberwithin{case}{thm}

\def\star(#1,#2){{\bf (\textasteriskcentered)$_{#1,#2}$}}

\def\eps{\epsilon}

\def\Fset{\mathcal{F}}
\def\HH{\mathcal{H}}
\def\K{\mathcal{K}}
\def\L{\mathcal{L}}
\def\Mset{\mathcal{M}}

\def\P{\mathcal{P}}
\def\Q{\mathcal{Q}}

\def\S{\mathcal{S}}
\def\T{\mathcal{T}}

\def\G{\mathcal{G}}

\def\N{\mathbb{N}}

\def\le{\leqslant}
\def\ge{\geqslant}
\def\ds{\displaystyle}
\def\ul{\underline}

\begin{document}

\title[Ordered graphs, ordered hypergraphs, and partitions]{Hereditary properties of partitions, ordered graphs and ordered hypergraphs}

\author{J\'ozsef Balogh}
\address{Department of Mathematics\\ University of Illinois\\ 1409 W. Green Street\\ Urbana, IL 61801} \email{jobal@math.uiuc.edu}

\author{B\'ela Bollob\'as}
\address{Department of Mathematical Sciences\\ The University of Memphis\\ Memphis, TN 38152\\ and\\ Trinity College\\ Cambridge CB2 1TQ\\ England} \email{bollobas@msci.memphis.edu}

\author{Robert Morris}
\address{Department of Mathematical Sciences\\ The University of Memphis\\ Memphis, TN 38152} \email{rdmorris@memphis.edu}
\thanks{The first author was supported during this research by OTKA grant T049398 and NSF grant DMS-0302804, the second by NSF grant ITR 0225610, and the third by a Van Vleet Memorial Doctoral Fellowship}

\begin{abstract}
In this paper we use the Klazar-Marcus-Tardos method (see \cite{MT}) to prove that if a hereditary property of partitions $\P$ has super-exponential speed, then for every $k$-permutation $\pi$, $\P$ contains the partition of $[2k]$ with parts $\{\{i,\pi(i) + k\} : i \in [k]\}$. We also prove a similar jump, from exponential to factorial, in the possible speeds of monotone properties of ordered graphs, and of hereditary properties of ordered graphs not containing large complete, or complete bipartite ordered graphs.

Our results generalize the Stanley-Wilf Conjecture on the number of $n$-permutations avoiding a fixed permutation, which was recently proved by the combined results of Klazar~\cite{Klaz} and Marcus and Tardos~\cite{MT}. Our main results follow from a generalization to ordered hypergraphs of the theorem of Marcus and Tardos.
\end{abstract}

\maketitle

\section{Introduction}\label{S:intro}

In this paper we shall prove that a jump from exponential to factorial speed occurs for properties of combinatorial structures of various types. We request the reader's patience while we make the various definitions necessary to state our results.

An \emph{ordered hypergraph} $\HH = (V,E,<)$ is a hypergraph -- a set of vertices $V$ and edges $E \subset \{A : A \subset V, |A| \ge 2\}$ -- together with a linear order $<$ on its vertices. Note that we do not allow edges to be repeated, and that we do not allow edges to consist of a single vertex. An ordered hypergraph $\K = (U,F,<)$ is an \emph{induced sub-hypergraph} of $\HH$ if $U \subset V$ (with the induced ordering), and $F = \{e \cap U : e \in E, |e| \ge 2\}$. $\K$ is a \emph{sub-hypergraph} of $\HH$ if $U \subset V$ (again with the induced ordering), and $F \subset \{e \cap U : e \in E, |e| \ge 2\}$. Finally, $\K$ is \emph{contained} in $\HH$ if there exists a sub-hypergraph $\L = (U,D,<)$ of $\HH$, with $|D| = |F| = t$, say, and $f_i \subset d_i$ for each $i \in [t]$ (where $D = \{d_1, \ldots, d_t\}$ and $F = \{f_1, \ldots, f_t\}$).

A collection of ordered hypergraphs is called a \emph{property} if it is closed under order-preserving isomorphisms of the vertex set. A property of ordered hypergraphs $\P$ is called \emph{hereditary} if it is closed under taking induced sub-hypergraphs; it is called \emph{monotone} if it is closed under taking sub-hypergraphs; and it is called \emph{strongly monotone} if it is closed under containment. Observe that any strongly monotone property is
monotone, and any monotone property is hereditary.

An \emph{ordered graph} is a graph together with a linear order $<$ on its vertices; equivalently, it is an ordered hypergraph in which each edge has size exactly 2. The definitions of hereditary and monotone properties are therefore inherited (note that in this case the definitions of monotone and strongly monotone coincide).

A \emph{partition} of the set $[n] = \{1, \ldots, n\}$ is an (unordered) collection of disjoint, non-empty sets $\{A_1, \ldots, A_t\}$ such that $A_1 \cup \ldots \cup A_t = [n]$. It is easy to see that a partition may be thought of as an ordered graph in which each component is a clique, or as an ordered hypergraph in which the edges are pairwise disjoint. Thus we obtain the definition of a hereditary property of partitions. Since we have come some distance from the original definition, we remark that if $P = \{A_1, \ldots, A_t\}$ is a partition of $[n]$, and $S$ is a subset of $[n]$ with elements $s_1 < \ldots < s_k$, then the \emph{sub-partition of $P$ induced by $S$} is the following partition of $[k]$. First let $\{B_1, \ldots, B_t\}$ satisfy $i \in B_j$ if and only if $s_i \in A_j$; then delete the empty classes. A property of partitions is hereditary if it is closed under taking sub-partitions.

Now, given a property $\P$ of ordered hypergraphs, we write $\P_n$ for the collection of distinct (i.e., non-isomorphic) ordered hypergraphs on $n$ vertices in $\P$, and call the function $n \mapsto |\P_n|$ the \emph{speed} (or unlabelled speed) of $\P$. An analogous definition can be made for other combinatorial structures (e.g., graphs, posets, permutations).

We are interested in the (surprising) phenomenon that for many such structures, only very `few' speeds are possible. More precisely, there often exists a family $\Fset$ of functions $f : \N \to \N$ and another function $F : \N \to \N$ with $F(n)$ {\em much} larger than $f(n)$ for every $f \in \Fset$, such that if for each $f \in \Fset$ the speed is infinitely often larger than $f(n)$, then it is also larger than $F(n)$ for every $n \in \N$. Putting it concisely: the speed {\em jumps} from $\Fset$ to $F$.

The study of the speeds of monotone properties of labelled graphs was introduced over forty years ago by Erd\H{o}s~\cite{E}, and continued by Erd\H{o}s, Kleitman and Rothschild~\cite{EKR}, Erd\H{o}s, Frankl and R\"odl~\cite{EFR}, Kolaitis, Pr\"omel and Rothschild~\cite{KPR}, Kleitman and Winston~\cite{KW}, Hundack, Pr\"omel and Steger~\cite{HPS} and more recently Balogh, Bollob\'as and Simonovits~\cite{BBS}. A new direction was initiated by Scheinerman and Zito~\cite{SZ}, who were the first to study hereditary properties of graphs with speeds below $n^n$. A little later, considerably stronger results were proved by Balogh, Bollob\'as and Weinreich~\cite{BBW}, \cite{BBW4}. In the range $|\P_n| = 2^{cn^2}$, the main results were proved by Alekseev~\cite{Alekseev}, Bollob\'as and Thomason~\cite{BTbox}, \cite{BT1}, and Pr\"omel and Steger~\cite{PS3}, \cite{PS4}, \cite{PS5}. For a review of the early results, see Bollob\'as~\cite{ICM}. Hereditary properties of other combinatorial structures have not yet been studied in such great detail, but it is likely that many more beautiful theorems await discovery.

In this paper we shall prove that a jump of this type, from exponential to factorial speed, occurs for strongly monotone properties of ordered hypergraphs. As a result of this, we shall be able to prove similar jumps for hereditary properties of partitions, monotone properties of ordered graphs, and hereditary properties of ordered graphs not containing arbitrarily large complete, or complete bipartite ordered graphs. As we shall see, each of these theorems is a generalization of the Stanley-Wilf Conjecture (Theorem A), proved recently by the combined results of Klazar~\cite{Klaz} and Marcus and Tardos~\cite{MT}.

Before we begin, we should remark that our main theorem has been proved independently (and at around the same time) by Klazar and Marcus~\cite{KM}. Although we were unaware of their work until after ours was completed, we should note also that many of the ideas in this paper were inspired by the earlier work of Klazar~\cite{Klaz} and of Marcus and Tardos~\cite{MT}. The reader may also wish to refer to some other papers of Klazar~\cite{Klaz2}, \cite{Klaz3} which we later discovered contain some of the ideas (though none of the results) below.

For each $n \in \N$, let $\Pi_n$ denote the collection of all permutations of $[n]$, and let $\Pi = \bigcup_n \Pi_n$. Also, if $\P$ is a property of ordered hypergraphs, and $k, \ell \in \N \cup \{\infty\}$, let $$\P^{(k,\ell)} = \{G \in \P : \Delta(G) \le k\textup{ and } |e| \le \ell \textup{ for every }e \in E(G)\}$$ denote the sub-property consisting of the ordered hypergraphs in which each vertex is contained in at most $k$ edges, and each edge has size at most $\ell$. Note that if $\P$ is hereditary, or monotone, or strongly monotone, then so is $\P^{(k,\ell)}$.

Finally, if $\pi \in \Pi_k$, let $H(\pi)$ denote the ordered hypergraph on vertex set $[2k]$ and with edge set $\{ \{i, \pi(i) + k\} : i \in [k]\}$. We shall also write $H(\pi)$ for the ordered graph with the same vertex and edge sets, and for the partition $\{\{i, \pi(i) + k\} : i \in [k]\}$ of $[2k]$. It will always be clear which of these $H(\pi)$ is.

We are now ready to state the main result of this paper, which was conjectured by Klazar in \cite{Klaz3}, and has been proved independently by Klazar and Marcus~\cite{KM}.

\begin{thm}\label{hypergraphs}
Let $\P$ be a strongly monotone property of ordered hypergraphs. If for every constant $c > 0$ there exists an $N = N(c) \in \N$ such that $|\P_N| > c^N$, then $\P$ contains the ordered hypergraph $H(\pi)$ for every $\pi \in \Pi$, and hence $$|\P_n| \: \ge \: |\P^{(1,2)}_n| \: \ge \: \ds\sum_{k=0}^{\lfloor n/2 \rfloor} {n \choose {2k}} k! \: = \: n^{n/2 + o(n)}$$ for every $n \in \N$. This lower bound is best possible, and there is a unique strongly monotone property of ordered hypergraphs with this speed.
\end{thm}

We remark that Scheinerman and Zito [25] proved that a similar jump, from exponential to factorial speed, exists for hereditary properties of labelled graphs. For more involved results, see \cite{BBW}.

Our proof of Theorem~\ref{hypergraphs} is based on the ideas of Klazar~\cite{Klaz}, \cite{Klaz3}, and of Marcus and Tardos~\cite{MT}. Theorem A, below, was proved in two stages: first Klazar~\cite{Klaz} showed that the theorem was a consequence of a conjecture of F\"uredi and Hajnal~\cite {FH}; then Marcus and Tardos~\cite{MT} proved that conjecture. In Theorem~\ref{genMT} we shall prove a generalization of the theorem of Marcus and Tardos (Theorem B). We shall then deduce Theorem~\ref{hypergraphs} using the method of Klazar (see \cite{Klaz} and \cite{Klaz3}).

We shall state Theorem~\ref{genMT} in terms of $(0,1)$--matrices, though it can equally be thought of as a theorem about ordered hypergraphs. To simplify the statement, we need a little notation.

Let $k \in \N$. If $A = (a_{ij})$ and $B = (b_{ij})$ are $k \times k$ matrices, we shall write $(A,B)$ for the $k \times 2k$ matrix obtained by putting $A$ in front of $B$. Thus $(A,B)_{ij}$ is $a_{ij}$ if $j \le k$ and $b_{i(j-k)}$ if $j \ge k+1$. Call two matrices $C$ and $D$ equivalent (and write $C \sim D$) if $D$ is obtained from $C$ by permuting its rows. Let $\Mset(k)$ denote the set of equivalence classes (with respect to $\sim$) in the family of all matrices of the form $(K,L)$, where $K$ and $L$ are $k \times k$ permutation matrices. Note that every such matrix $(K,L)$ is equivalent to a unique matrix $(I,M)$, where $I = (\delta_{ij})$ is the $k \times k$ identity matrix, so $|\Mset(k)| = k!$.

Finally, if $P$ and $Q$ are $(0,1)$--matrices, then we say that $Q$ is a \emph{sub-matrix} of $P$ if $Q$ is obtained from $P$ by deleting rows and columns. We say that $P$ \emph{contains} $Q = (q_{ij})$ if there exists a sub-matrix $R = (r_{ij})$ of $P$, the same size as $Q$, with $r_{ij} = 1$ whenever $q_{ij} = 1$. If we associate an ordered hypergraph $\HH$ with a $(0,1)$--matrix whose rows are the indicator functions of the edges of $\HH$, and consider two matrices to be the same if they are equivalent, then this concept coincides with hypergraph containment defined above. If $P$ does not contain $Q$ then we say that $P$ \emph{avoids} $Q$.

The following theorem is a reformulation of Klazar's conjecture C4 in \cite{Klaz3}. It has also been proved independently by Klazar and Marcus~\cite{KM}.

\begin{thm}\label{genMT}
Let $k \in \N$. There exists a constant, $c_k$, such that if $m,n \in \N$ and $A$ is an $m \times n$ $(0,1)$--matrix satisfying
\begin{enumerate}
\item[$(i)$] at least $c_kn$ of the entries of $A$ are $1$, and
\item[$(ii)$] each of the rows of $A$ are different,
\end{enumerate}
then $A$ contains some member of each class of $\Mset(k)$.
\end{thm}

\noindent Notice that Theorem~\ref{genMT} still holds if condition $(ii)$ is replaced by the
condition\\[+1ex]
\indent $(ii')$ each of the rows of $A$ has at least $2k$ of its entries $1$,\\[+1ex]
since if $A$ satisfies $(ii')$, and any row occurs at least $k$ times in $A$, then $A$ contains every $k \times 2k$ $(0,1)$-- matrix.

\begin{rmk}
F\"uredi and Hajnal~\cite {FH} proved that the extremal number of $1$'s possible in an $n \times n$ $(0,1)$--matrix avoiding the matrix $$S_1 = \left( \begin{array}{cccc} 1 & 0 & 1 & 0 \\ 0 & 1 & 0 & 1 \end{array} \right)$$ is (up to a constant) $n\alpha(n)$, where $\alpha(n) \to \infty$ extremely slowly. A simple corollary of Theorem~\ref{genMT} (with $k = 2$) is that the extremal number of $1$'s if we avoid both $S_1$ and $$S_2 = \left( \begin{array}{cccc} 0 & 1 & 0 & 1 \\ 1 & 0 & 1 & 0 \end{array} \right)$$ is $O(n)$. For many more results along these lines, see Tardos~\cite{T}.
\end{rmk}

We shall now note two important (and immediate) consequences of Theorem~\ref{hypergraphs}. The first of them was conjectured by Klazar in \cite{Klaz2}, and the second was originally proved (although not stated!) by Klazar in \cite{Klaz} as a consequence on the F\"uredi--Hajnal Conjecture.

\begin{thm}\label{parts}
Let $\P$ be a hereditary property of partitions. If for every constant $c > 0$ there exists an $N = N(c) \in \N$ such that $|\P_N| > c^N$, then $\P$ contains the partition $H(\pi)$ for every $\pi \in \Pi$, and hence
$$|\P_n| \: \ge \: |\P^{(1,2)}_n| \: \ge \: \ds\sum_{k=0}^{\lfloor n/2 \rfloor} {n \choose {2k}} k! \: = \: n^{n/2 + o(n)}$$ for every $n \in \N$. This lower bound is best possible, and there is a unique hereditary property of partitions with this speed.
\end{thm}

\begin{thm}\label{mono}
Let $\P$ be a monotone property of ordered graphs. If for every constant $c > 0$ there exists an $N = N(c) \in \N$ such that $|\P_N| > c^N$, then $\P$ contains the ordered graph $H(\pi)$ for every $\pi \in \Pi$, and hence
$$|\P_n| \: \ge \: |\P^{(1,2)}_n| \: \ge \: \ds\sum_{k=0}^{\lfloor n/2 \rfloor} {n \choose {2k}} k! \: = \: n^{n/2 + o(n)}$$ for every $n \in \N$. This lower bound is best possible, and there is a unique monotone property of ordered graphs with this speed.
\end{thm}

For $t \in \N$, let $K_t$ denote the complete ordered graph on $t$ vertices, and let $K_{t,t}$ denote the complete ordered bipartite graph on $[2t]$ with edge set $E(K_{t,t}) = \{\{i,j\} : i \le t < j\}$. We shall deduce the following theorem from Theorem~\ref{mono}.

\begin{thm}\label{noKt}
Let $\P$ be a hereditary property of ordered graphs such that for some $t \in \N$, neither $K_t$ nor $K_{t,t}$ is in $\P$. If for every constant $c > 0$ there exists an $N = N(c) \in \N$ such that $|\P_N| > c^N$, then $\P$ contains the ordered graph $H(\pi)$ for every $\pi \in \Pi$, and hence
$$|\P_n| \: \ge \: |\P^{(1,2)}_n| \: \ge \: \ds\sum_{k=0}^{\lfloor n/2 \rfloor} {n \choose {2k}} k! \: = \: n^{n/2 + o(n)}$$
for every $n \in \N$. This lower bound is best possible, and there is a unique hereditary property containing neither $K_t$ nor $K_{t,t}$ with this speed.
\end{thm}

We conjecture that Theorems~\ref{hypergraphs}, \ref{parts}, \ref{mono} and \ref{noKt} have the following common generalization.

\begin{conj}\label{hyperconj}
Let $\P$ be a hereditary property of ordered hypergraphs. If for every constant $c > 0$ there exists an $N = N(c) \in \N$ such that $|\P_N| > c^N$, then
$$|\P_n| \: \ge \: \ds\sum_{k=0}^{\lfloor n/2 \rfloor} {n \choose {2k}} k! \: = \: n^{n/2 + o(n)}$$ for every $n \in \N$.
\end{conj}

The following statement is a special case of Conjecture~\ref{hyperconj}, but still generalizes Theorems~\ref{parts}, \ref{mono} and \ref{noKt} (since partitions can be represented by ordered graphs whose components are complete graphs), and would be very interesting in its own right. It was in fact our main motivation for studying ordered hypergraphs and partitions.

\begin{conj}\label{orderconj}
Let $\P$ be a hereditary property of ordered graphs. If for every constant $c > 0$ there exists an $N = N(c) \in \N$ such that $|\P_N| > c^N$, then $$|\P_n| \: \ge \: \ds\sum_{k=0}^{\lfloor n/2 \rfloor} {n \choose {2k}} k! \: = \: n^{n/2 + o(n)}$$ for every $n \in \N$.
\end{conj}

\begin{rmk}
Note that in both conjectures the lower bounds, if true, are best possible (by Lemma~\ref{bound}, below). However it is \emph{not} true that, under the conditions of the conjectures, $\P$ must contain the ordered graph $H(\pi)$ on $[2k]$ with edge set $\{\{i,\pi(i)+k\} : i \in [k]\}$ for every $k \in \N$ and $\pi \in \Pi_k$. To see this, call an ordered graph $G$ on $[n]$ a \emph{co-matching} if $\{x_1,y_1\},\{x_2,y_2\} \in {{[n]} \choose 2} \setminus E(G)$ implies that $|\{x_1,y_1\} \cap \{x_2,y_2\}| = 0$ or $2$, and call $G$ a \emph{star-matching} if (say) $\{x_1,y_1\} \in E(G)$ implies $\{x_1,y_2\} \in E(G)$ for every $y_1 \le y_2 \le n$. The collection of all co-matchings and the collection of all star-matchings are hereditary properties of ordered graphs with super-exponential speeds, but neither contains all the graphs $H(\pi)$.
\end{rmk}

For further details on the possible speeds of hereditary properties of ordered graphs see \cite{BBMord}, which considers such properties with speed below $2^n$, and also those with speed above $2^{\eps n^2}$.

The rest of the paper is organised as follows. In Section~\ref{KMT} we shall state the Klazar-Marcus-Tardos and Marcus-Tardos Theorems, and show that the former is implied by each of Theorems~\ref{parts}, \ref{mono} and \ref{noKt}; in Section~\ref{genMTsec} we shall prove Theorem~\ref{genMT} using the Marcus-Tardos theorem; in Section~\ref{hypersec} we shall deduce Theorem~\ref{hypergraphs} from Theorem~\ref{genMT}; and in Section~\ref{rest} we shall deduce Theorems~\ref{parts}, \ref{mono} and \ref{noKt}.

\section{The Klazar-Marcus-Tardos and Marcus-Tardos theorems}\label{KMT}

We begin by recalling the theorems of Marcus and Tardos~\cite{MT}.

Given $n \in \N$, we shall call a permutation of $[n]$ an $n$-permutation. An $n$-permutation $\pi$ is said to \emph{contain} a $k$-permutation $\sigma$ if there are integers $1 \le a(1) < \ldots < a(k) \le n$ such that $\pi(a(i)) < \pi(a(j))$ if and only if $\sigma(i) < \sigma(j)$. Otherwise $\pi$ is said to \emph{avoid} $\sigma$. A property of permutations is a collection of permutations, closed under isomorphism. A property of permutations is said to be hereditary if it is also closed under containment.

The following theorem was conjectured by Stanley and Wilf around 1992 (see \cite{Arratia}, \cite{Bona}, \cite{MT}), and proved by Marcus and Tardos in 2004 (Corollary 2 of \cite{MT}), using a theorem of Klazar~\cite{Klaz}. This result is usually known as the Stanley-Wilf Conjecture, but we shall refer to it as the Klazar-Marcus-Tardos Theorem, or simply as Theorem A.

\begin{thma}
Let $\P$ be a hereditary property of permutations. Either $\P$ is the set $\Pi$ of all permutations, so $|\P_n| = n!$ for every $n \in \N$, or there exists a constant $c = c(\P)$ such that $|\P_n| \le c^n$ for every $n \in \N$.
\end{thma}

We also state here the theorem of Marcus and Tardos (Theorem 1 of \cite{MT}), which was originally conjectured by F\"uredi and Hajnal in \cite{FH}. In Section~\ref{genMTsec} we shall use it to prove Theorem~\ref{genMT}.

\begin{thmb}
For every permutation matrix $M$, there exists a constant $C = C(M)$ such that any $n \times n$ $(0,1)$--matrix with at least $Cn$ of its entries $1$ contains $M$.
\end{thmb}

In this section we shall show that the simplest case of Conjectures~\ref{hyperconj} and \ref{orderconj} (the case in which every $G \in \P$ is an ordered graph with maximum degree at most one) is equivalent to the Klazar-Marcus-Tardos Theorem, and deduce that our main results generalize that theorem. We start however by proving the following lemma, which gives the final implication of Theorems~\ref{hypergraphs}, \ref{parts}, \ref{mono} and \ref{noKt}.

\begin{lemma}\label{bound}
Let $\P$ be a hereditary property of ordered hypergraphs. If $H(\pi) \in \P$ for every $\pi \in \Pi$, then $$|\P_n| \ge |\P^{(1,2)}_n| \ge \ds\sum_{k=0}^{\lfloor n/2 \rfloor} {n \choose {2k}} k!$$ for every $n \in \N$. Moreover, there is a unique hereditary property containing every $H(\pi)$ with this speed.
\end{lemma}

\begin{proof}
Given integers $k,n \in \N$, a subset $A \subset [n]$ of size $2k$ (with elements $a(1) < \ldots < a(2k)$ say), and permutation $\pi \in \Pi_k$, define $G(n,A,\pi)$ to be the ordered hypergraph on vertex set $[n]$, and with edge set $\{\{a(i), a(\pi(i)+k)\} : i \in [k]\}$. Let $\P$ be a hereditary property of ordered hypergraphs with $H(\pi) \in \P$ for every $\pi \in \Pi$. We shall show that $G(n,A,\pi) \in \P$ for every such $n$, $A$ and $\pi$.

Indeed, let $n$, $A$ and $\pi$ be as described, let $X$ be the set of isolated vertices in $G = G(n,A,\pi)$, let $Y = \{v \in X : v < a_k\}$ and let $Z = X \setminus Y$. Suppose $|Y| = r$ and $|Z| = s$, and consider an ordered graph $H$ on $[-s+1,n+r]$ formed by adding to $G$ an arbitrary matching between the vertices $\{-s+1, \ldots, 0\}$ and $Z$, and an arbitrary matching between the vertices $\{n+1, \ldots, n+r\}$ and $X$. It is easy to see that $H$ is isomorphic to $H(\sigma)$ for some $\sigma \in \Pi_{k+r+s}$, so $H \in \P$, and that $G$ is an induced subgraph of $H$, so $G \in \P$.

Thus $\P$ contains the ordered hypergraph $G(n,A,\pi)$ for every $k \in \N$, $\pi \in \Pi_k$, $n \ge 2k$ and $A \subset [n]$ with $|A| = 2k$, and hence
\begin{eqnarray*}
|\P_n| \; \ge \; |\P^{(1,2)}_n| & \ge & |\{(k,\pi,A) : 2k \le n\textup{, }\pi \in \Pi_k\textup{ and }A \subset [n]^{(2k)}\}|\\
& = & \ds\sum_{k=0}^{\lfloor n/2 \rfloor} {n \choose {2k}}
k!.\end{eqnarray*}

Finally, note that the collection $\Q = \{G(n,\pi,A) : \pi \in \Pi_k$ and $A \subset [n]^{(2k)}$ for some $k \le n/2\}$ forms a hereditary property of ordered hypergraphs, and $H(\pi) \in \Q$ for every $\pi \in \Pi$. By the argument above, if $H(\pi) \in \P$ for every $\pi \in \Pi$ then $\Q \subset \P$, so $\Q$ the unique such hereditary property of ordered hypergraphs with this speed.
\end{proof}

We shall now show that the simplest case of Conjecture~\ref{hyperconj} follows from the Klazar-Marcus-Tardos Theorem. We shall not need this to prove our main results, but the proof is short and has some independent value. Here we use $\G$ to denote a property of ordered graphs, to distinguish it from a property $\P$ of permutations.

\begin{thm}\label{deg1}
Let $\G$ be a hereditary property of ordered graphs of maximal degree at most $1$. If for every constant $c > 0$ there exists an $N = N(c) \in \N$ such that $|\G_N| > c^N$, then $\G$ contains the ordered graph $H(\pi)$ for every $\pi \in \Pi_k$, and hence $$|\G_n| \ge |\G^{(1,2)}_n| \ge \ds\sum_{k=0}^{\lfloor n/2 \rfloor} {n \choose {2k}} k!$$ for every $n \in \N$.
\end{thm}

\begin{proof}
Let $\G$ be a hereditary property of ordered graphs of maximal degree at most one, and suppose that for every constant $c_1 > 0$ there exists an $N = N(c_1) \in \N$ such that $|\G_N| > c_1^N$. Given an ordered graph $G \in \G$, we define a permutation $\phi(G)$. Suppose $G$ has $k$ edges, $e_1 = \{a_1,b_1\}, \ldots, e_k = \{a_k,b_k\}$, (where $a_i < b_i$ for each $i \in [k]$), ordered by their left-endpoints, i.e., $a_i < a_j$ if and only if $i < j$ (recall that $\Delta(G) \le 1$). Let $\pi$ be the $k$-permutation such that $b_{\pi(1)} < \ldots < b_{\pi(k)}$, and define $\phi(G) = \pi^{-1}$.

Let $\P = \{\phi(G) : G \in \G\}$. Since $\G$ is hereditary, so is $\P$, since removing a vertex from a permutation corresponds to removing one of the endpoints of the corresponding edge. By the Klazar-Marcus-Tardos Theorem, either $\P = \Pi$ or there exists a constant $c$ such that $|\P_n| \le c^n$ for every $n \in \N$.

Suppose the latter, so $|\P_n| \le c^n$ for every $n \in \N$. Assuming $c > 2$, we claim that in this case $$|\G_n| \: \le \: \ds\sum_{k=0}^{\lfloor n/2 \rfloor} {n \choose k}{{n-k} \choose k}c^k \: < \: n {n \choose {\lfloor n/2 \rfloor}}^2 c^{\lfloor n/2 \rfloor} \: < \: (4\sqrt{c})^n.$$ To see this, simply note that any ordered graph $G$ of maximal degree at most one is determined by its order, its left-endpoint set, its right-endpoint set, and the permutation $\phi(G)$. Hence, setting $c_1 = 4\sqrt{c}$, we have a contradiction to our assumption that $|\G_N| > c_1^N$ for some $N \in \N$.

Next, suppose that $\P = \Pi$. We want to show that $\G$ contains $H(\pi)$ for every $\pi \in \Pi$, so let us fix $k \in \N$ and $\pi \in \Pi_k$. Let $\pi'$ be the $(k+1)$-permutation defined as follows: $\pi'(i) = \pi(i) + 1$ for each $i \in [k]$, and $\pi'(k+1) = 1$. Since $\P = \Pi$, we have $\pi' \in \P$, so for some $G \in \G$ we have $\pi' = \phi(G)$.

Now notice that in $G$, all left-endpoints occur to the left of \emph{all} right-endpoints, since $\pi'(k+1) = 1$. Therefore, letting $G'$ be subgraph of $G$ induced by the first $k$ left-endpoints and last $k$ right-endpoints, we have $G' = H(\pi)$.
\end{proof}

\begin{rmk}
In fact, for any permutation $\pi \in \Pi_k$, the number of ordered graphs $G$ of order $n$ and maximal degree at most one with $\phi(G) = \pi$ is at most ${n \choose k}Cat(k)$, where $Cat(k) = \frac{1}{k+1}{{2k} \choose k}$ is the $k^{th}$ Catalan number. To see this, we use the fact that there are exactly $Cat(k)$ legal sequences of $k$ left- and $k$ right-brackets (i.e., in any initial segment of the sequence there are at least as many left-brackets as right-brackets). Given any ordered graph $G \in \G$, we can define a corresponding sequence of brackets, $\psi(G)$, by taking a left-bracket for every vertex which is the left-endpoint of an edge, and a right-bracket for every right-endpoint.

Now, given $n$, $\phi = \phi(G)$, $\psi = \psi(G)$ and the (even-sized) subset $A = \{v \in [n] : d_G(v) = 1\}$, it is simple to reconstruct $G$: if the elements of $A$ are $a(1) < \ldots < a(2k)$, the left brackets of $\psi$ lie in positions $1 \le s(1) < \ldots < s(k) \le 2k$ and the right brackets lie in positions $1 \le t(1) < \ldots < t(k) \le 2k$ (so $\{s(1), \ldots, s(k), t(1), \ldots, t(k)\} = [2k]$), then the edge set is $\{\{a(s(i)),a(t(\phi(i) ) ) \} : i \in [k]\}$. Note that although for many permutation - bracket sequence pairs $(\phi, \psi)$ no ordered graph $G$ has $\phi(G) = \phi$ and $\psi(G) = \psi$ (for example, $\phi = 21$ and $\psi =$ [( ) ( )]), for the identity permutation all $Cat(k)$ bracket pairs are realised.
\end{rmk}

We proved Theorem~\ref{deg1} using Theorem A; we now prove the reverse implication. It will follow almost immediately that Theorem $i$ implies Theorem A for $i =$~\ref{parts}, \ref{mono} and \ref{noKt}.

\begin{lemma}\label{genA}
Theorem~\ref{deg1} implies Theorem A.
\end{lemma}

\begin{proof}
Let $\P$ be a non-trivial hereditary property of permutations (i.e., different from $\Pi$), and assume that Theorem~\ref{deg1} holds. Let the ordered graphs $G(n,A,\pi)$ be as defined above, and let
$$\G = \{G(n,A,\pi) : n \in \N, \: A \subset [n], \: \pi \in \P, \: |A| = 2|\pi|\}.$$
Because $\P$ is hereditary, $\G$ is also hereditary, since removing an isolated vertex from $G(n,A,\pi)$ gives $G(n-1,A',\pi)$ (for some $A' \subset [n-1]$), and removing a non-isolated vertex corresponds to removing an element from $\pi$.

Since $\P \neq \Pi$ there exists some $\pi \notin \P$, and by definition $\G$ does not contain $H(\pi)$. Hence, by Theorem~\ref{deg1}, there exists a constant $c > 0$ such that $|\G_n| \le c^n$ for every $n \in \N$. But $|\G_{2n}| \ge |\P_n|$ for every $n \in \N$, so $|\P_n| \le c^{2n}$ for every $n \in \N$.
\end{proof}

We can now deduce that our main theorems do indeed generalize the Klazar-Marcus-Tardos Theorem.

\begin{cor}
Each of the Theorems~\ref{parts}, \ref{mono} and \ref{noKt} implies Theorem A.
\end{cor}

\begin{proof}
To show that Theorem~\ref{parts} and Theorem~\ref{noKt} imply Theorem A, it suffices to observe that any hereditary property of ordered graphs of maximal degree at most one may be viewed as a hereditary property of partitions (with part sizes at most 2), or as a hereditary property of ordered graphs containing no $K_3$ and no $K_{2,2}$. The result then follows by Lemma~\ref{genA}.

To show that Theorem~\ref{mono} implies Theorem A, let $\P$ be a hereditary property of ordered graphs of maximal degree one, and consider the minimal monotone property of ordered graphs $\P'$ containing $\P$. If $\P'$ contains the ordered graph $H(\pi)$ (for some $\pi \in \Pi$) then so does $\P$. Otherwise $|\P'_n| \le c^n$ for some $c > 0$ and every $n \in \N$ by Theorem~\ref{mono}, and hence $|\P_n| \le |\P'_n| \le c^n$ for every $n \in \N$. The result again follows by Lemma~\ref{genA}.
\end{proof}

\section{Proof of the generalized Marcus-Tardos Theorem}\label{genMTsec}

In this section we shall prove Theorem~\ref{genMT}. Recall that by Theorem B, for each permutation matrix $M$, there exists a constant $C(M)$ such that any $n \times n$ $(0,1)$--matrix with at least $C(M)n$ of its entries $1$ contains $M$. For each $k \in \N$, let $C(k)$ be the constant obtained in the Theorem B for $k \times k$ matrices, i.e., $C(k) = \operatorname{max}\{C(M) : M$ a $k \times k$ permutation matrix$\}$. We shall give our bounds on $c_k$ in terms of $C(k)$.

To obtain Theorem~\ref{genMT}, we use Theorem B to prove it in the case that the rows each have a bounded number of $1$'s, and then use this result and the method of Marcus and Tardos~\cite{MT} to prove the general case. First however, we need to show that Theorem B implies Theorem~\ref{genMT} in the case that each row has exactly two $1$'s; in fact these statements are equivalent.

\begin{lemma}\label{MTdeg2}
Let $f: \N \to \N$ be any function. The following statements satisfy $(i) \Rightarrow (ii) \Rightarrow (iii)$.\\
(i) For each $k,m,n \in \N$, any $m \times n$ $(0,1)$--matrix with at least $f(k)n$ of its entries $1$, and at most two $1$'s in each row, and with each row different, contains a member of each class of $\Mset(k)$.\\
(ii) Theorem B holds with $C(M_k) = f(k)$ for each $k \times k$ permutation matrix $M_k$.\\
(iii) For each $k,m,n \in \N$, any $m \times n$ $(0,1)$--matrix with at least $(2f(k+1)+1)n$ of its entries $1$, and at most two $1$'s in each row, and with each row different, contains a member of each class of $\Mset(k)$.
\end{lemma}

\begin{proof}
First we shall prove that $(i)$ implies $(ii)$. Let $k,m,n \in \N$, $M = (m_{ij})$ be a $k \times k$ permutation matrix, and $A = (a_{ij})$ be an $n \times n$ $(0,1)$--matrix with at least $f(k)n$ $1$'s, with $f(k)$ given by $(i)$. We wish to show that $A$ contains $M$. Suppose without loss of generality that there are more $1$'s above the top-left/bottom-right diagonal than below it (otherwise replace $A$ and $M$ by $A^T$ and $M^T$). Let the number of pairs $(i,j)$ for which $i \le j$ and $a_{ij} = 1$ be $m$, and label them $e_1, \ldots, e_m$ arbitrarily. We define an $m \times n$ $(0,1)$--matrix $B = (b_{ij})$ with at most two $1$'s in each row, and with each row different, by letting $b_{ij} = 1$ if and only if vertex $j$ is an endpoint of $e_i$. Note that at least $f(k)n$ of the entries of $B$ are $1$.

Applying $(i)$ to $B$, we see that $B$ must contain a matrix $(K,L) \sim (I,M)$, where $(I,M)$ is the $k \times 2k$ matrix obtained by putting the identity matrix in front of $M$ (so $(I,M)_{ij} = \delta_{ij}$ and $(I,M)_{i(j+k)} = m_{ij}$ for each $j \in [k]$). Suppose $(K,L)$ occurs in columns $a_1 < \ldots < a_k < b_1 < \ldots < b_k$ of $B$. Then $M$ occurs in the intersection of the rows $a_1,\dots,a_k$ and the columns $b_1, \ldots, b_k$ of $A$, and so we are done.

The proof that $(ii)$ implies $(iii)$ is similar. Again let $k,m,n \in \N$, $M = (m_{ij})$ be a $k \times k$ permutation matrix, and let $B = (b_{ij})$ be an $m \times n$ $(0,1)$--matrix with at most two $1$'s in each row, each row different, and at least $(2f(k+1)+1)n$ of its entries $1$. It will suffice to show that $B$ contains some matrix $(K,L) \sim (I,M)$.

We produce from $B$ an $n \times n$ $(0,1)$--matrix $A = (a_{ij})$, by letting $a_{ij} = 1$ if and only if $i < j$ and $b_{ri} = b_{rj} = 1$ for some $r \in [m]$. Note that at least $f(k+1)n$ of the entries of $A$ are $1$. Applying $(ii)$ to $A$, we see that $A$ must contain the ($k+1) \times (k+1)$ matrix $M'$, formed by putting $M$ in the top right-hand corner, and a single $1$ in the bottom left-hand corner. Thus $M'_{i(j+1)} = m_{ij}$ for $i,j \in [k]$, $M'_{(k+1)1} = 1$, and $M'_{ij} = 0$ otherwise. Suppose $M'$ occurs in rows $a_1 < \ldots < a_{k+1}$ and columns $b_1 < \ldots < b_{k+1}$ of $A$, and that the $1$'s corresponding to $1$-entries of $M$ correspond to rows $r_1 < \ldots < r_k$ of $B$. Since $A_{a_{k+1}b_1} = 1$ and $A$ is upper triangular, $a_{k+1} < b_1$. Therefore some $(K,L) \sim (I,M)$ occurs in the intersection of the rows $r_1,\dots,r_k$ and the columns $a_1, \ldots, a_k, b_2, \ldots, b_{k+1}$ of $B$, and we are again done.
\end{proof}

Since Theorem B holds with $C(M) = C(k)$, we have the following immediate corollary.

\begin{cor}\label{MTdeg2(c)}
Let $k \in \N$. Any $m \times n$ $(0,1)$--matrix with at least $(2C(k+1)+1)n$ of its entries $1$, and at most two $1$'s in each row, and with each row different, contains a member of each class of $\Mset(k)$.
\end{cor}

To prove the case where the rows have a bounded number of $1$'s, we shall use the following trivial observation.

\begin{lemma}\label{match}
Let $G$ be a bipartite graph with parts $A$ and $B$. Suppose $d(v) \ge 1$ for each $v \in A$, and $d(v) \le m$ for each $v \in B$. Then there exists a matching in $G$ of size at least $|A|/m$.
\end{lemma}

\begin{proof}
If $d(v) > 1$ for any $v \in A$ then remove all but one of the edges joined to $v$. We now have a family of stars, each centred in $B$ and of order at most $m+1$. Take one edge from each.
\end{proof}

For each $D \in \N$, define $g_D : \N \to \N$ by $g_D(x) = \sum_{i=0}^{D-1} (i+1){x \choose i}$. Note that $g_D(x) < 2D{x \choose {D-1}}$ if $x > 3D$. Let $\HH$ be an ordered hypergraph on $[n]$ in which every edge has size at most $D$. For each vertex $v \in [n]$, the $2$-degree of $v$ in $\HH$ is $d^{(2)}_{\HH}(v) = |\{u \in [n] : \{u,v\} \subset E$ for some $E \in E(\HH)\}|$. Suppose a vertex of $2$-degree $x$ is removed from $\HH$; by how much can $\|\HH\| = \sum_{E \in E(\HH)} |E|$ decrease? For each $0 \le i \le D-1$, $v$ is contained in at most ${x \choose i}$ edges $E$ of $\HH$ of size $i+1$, and each of these must be removed entirely if $E \setminus v$ is also an edge. Hence the maximum possible decrease in $\|\HH\|$ is $g_D(x)$.

For each $k \in \N$ let $C_1(k) = 2C(k+1) + 1$.

\begin{lemma}\label{MTdegbdd}
Let $k,m,n,D \in \N$, with $D \ge 2$, and let $g_D: \N \to \N$ be as defined above. Any $m \times n$ $(0,1)$--matrix with at least $g_D\left( D(D - 1)C_1(k) \right)n$ of its entries $1$, and at most $D$ of the entries of each row $1$, and with each row different, contains a member of each class of $\Mset(k)$.
\end{lemma}

\begin{proof}
We shall use Corollary~\ref{MTdeg2(c)}. Let $k,m,n,D \in \N$, $M = (m_{ij})$ be a $k \times k$ permutation matrix, and $A = (a_{ij})$ be an $m \times n$ $(0,1)$--matrix with at least $g_D\left(D(D - 1) C_1(k) \right)n$ of its entries $1$, at most $D$ of the entries in each row $1$, and each row different. We shall show that $A$ contains a matrix in the equivalence class of $(I,M)$.

Consider the ordered hypergraph $\HH$ on vertex set $[n]$ with edge set $\{E_i : i \in [m]$ and $a_{ij} = 1 \Leftrightarrow j \in E_i\}$, so the rows of $A$ are the indicator functions of the edges. Note that $\|\HH\| = \sum_i |E_i| \ge g_D\left( D(D - 1) C_1(k) \right)n$. We first wish to find a subset $S$ of $[n]$ in which there are at least ${D \choose 2}C_1(k)|S|$ distinct pairs $\{i,j\}$, each contained in some edge of $\HH$. If $d^{(2)}_{\HH}(v) < D(D - 1)C_1(k)$ for some vertex $v$ then, by the comments above, removing $v$ from the ordered hypergraph causes $\|\HH\|$ to decrease by at most $g_D\left( D(D - 1)C_1(k) - 1 \right) < g_D\left( D(D - 1)C_1(k) \right)$. Thus removing $v$ causes the density of edges in the ordered hypergraph to increase.

Thus, if we repeatedly remove vertices of minimal $2$-degree from $\HH$, we must eventually produce an ordered hypergraph $\HH'$ on vertex set $S$ in which every vertex $v$ has $d^{(2)}_{\HH'}(v) \ge D(D - 1)C_1(k)$. By counting degrees, there are at least ${D \choose 2}C_1(k)|S|$ distinct pairs $\{i,j\} \subset S$, each contained in some edge of $\HH'$, and hence of $\HH$.

Let $\T = \{\{i,j\} \subset S : i,j \in E$ for some $E \in \HH\}$ be the set of such pairs. Let $B$ be the bipartite graph on sets $\T$ and $E(\HH)$ with edges corresponding to containment (i.e., $(\{i,j\},E)$ is an edge of $B$ iff $\{i,j\} \subset E$). By Lemma~\ref{match}, there exists a matching $W$ in $B$ of size $t$, with $t \ge |\T|/{D \choose 2} \ge C_1(k)|S|$, since each edge has order at most $D$. Let $\T_1$ be the set of endpoints of $W$ lying in $\T$.

Let $s = |S|$, and define $P$ to be the $t \times s$ $(0,1)$--matrix in which the columns correspond to elements of $S$, and the rows are the indicator functions of the edges in $\T_1$. In $P$ all rows have exactly two $1$'s, all rows are different, and  $2t \ge 2C_1(k)s$ of its entries are $1$, so by Corollary~\ref{MTdeg2(c)}, $P$ contains some matrix $(K,L)$ equivalent to $(I,M)$.

Fix a copy of $(K,L)$ in $P$, and let the columns of $P$ containing this $(K,L)$ be those corresponding to the vertices $b(1) < \ldots < b(2k)$ of $\HH$ (note that here $b(i)$ is the \emph{original} labelling of the vertex in $[n]$). We claim that the corresponding columns of $A$ contain a matrix equivalent to $(K,L)$. To see this, let the rows of $P$ containing the same copy of $(K,L)$ be $a(1) < \ldots < a(k)$, and let $p(i)$ be the pair in $\T_1$ corresponding to row $a(i)$ for $1 \le i \le k$. For each of these pairs $p(i)$, choose the edge $e(i) \in E(\HH)$ it was matched to by $W$. We have thus found $k$ distinct edges $e(1), \ldots, e(k)$, for which $p(i) \subset e(i)$. It follows immediately that the columns $b(1), \ldots, b(2k)$ of $A$ contain some matrix $(K',L') \sim (K,L) \sim (I,M)$. This completes the proof.
\end{proof}

We are now ready to prove Theorem~\ref{genMT}.

\begin{proof}[Proof of Theorem~\ref{genMT}]
For each $n,k \in \N$, let $f(n,k)$ be the largest number of $1$'s possible in an $m \times n$ $(0,1)$--matrix $A$ (where $m \in \N$ is arbitrary), with each row different, not containing any member of some class of $\Mset(k)$. We wish to show that $f(n,k) = O(n)$, where $k$ is fixed and $n \to \infty$. The proof (that $f(n,k) < c_k n$, where $c_k$ will be determined later) will be by induction on $n$. Note that since each row of $A$ is different, at most $n2^n$ of the entries of $A$ can be $1$'s. We shall choose $c_k > 2^{8k^3}$, so the statement $f(n,k) < c_kn$ is (trivially) true for $n \le 8k^3$.

Let $k,m,n \in \N$ with $n \ge 8k^3$, $M = (m_{ij})$ be a $k \times k$ permutation matrix, and $A$ be an $m \times n$ $(0,1)$--matrix, with each row different, not containing any matrix equivalent to $(I,M)$. Following the method of Marcus and Tardos, we want to divide $A$ up into `fat' and `skinny' blocks of size $1 \times t$ for some $t$. In preparation for this, we must remove the rows with few $1$'s. Let $t \in \N$ and let $D = (2k-1)t$. (We shall eventually set $t = 2k^2$, but we postpone choosing this value until it is clear why the choice is being made. Our argument up to that point works for any $t \in \N$.) By Lemma~\ref{MTdegbdd}, there are at most $g_D\left(D(D - 1)C_1(k) \right)n$ $1$'s in rows with at most $D$ of their entries $1$, otherwise $A$ would contain some $(K,L) \sim (I,M)$, contradicting our assumption. Let $A' = (a_{ij}')$ be the $m' \times n$ matrix obtained from $A$ by deleting the rows with at most $D$ entries $1$.

Now, let $q$ and $r$ satisfy $n = qt + r$, with $q \in \N$ and $r \in [t] $, and partition $A'$ into $qm'$ blocks of size $1 \times t$ and $m'$ blocks of size $1 \times r$ as follows. Let the $1 \times t$ blocks be $S_{ij} = (a_{i\ell}' : \ell \in [(j-1)t+1,jt])$, for each $i \in [m']$ and $j \in [q]$, and the $1 \times r$ blocks be $S_{i(q+1)} = (a_{i\ell}' : \ell \in [qt+1,n])$, for each $i \in [m']$. Define $B = (b_{ij})$ to be the $m' \times (q+1)$ $(0,1)$--matrix obtained by assigning the value $1$ to a block if any entry of the block is $1$. Thus, for $j \in [q]$, $b_{ij} = 0$ if and only if $a_{i\ell}' = 0$ for every $\ell \in [(j-1)t+1,jt]$, and similarly for $j = q+1$.\\[-1ex]

\noindent\ul{Claim 1}: $B$ contains no matrix equivalent to $(I,M)$.
\begin{proof}
This is Lemma 4 of Marcus and Tardos~\cite{MT}. To spell it out, assume $B$ contains such a matrix $P$, and for each $1$ it contains, choose an arbitrary non-zero entry from the corresponding blocks of $A'$. They represent a copy of $P$ in $A$, a contradiction.
\end{proof}

\noindent Call a block `fat' if at least $2k$ of its entries are $1$.\\[-1ex]

\noindent\ul{Claim 2}: There are at most ${t \choose {2k}}(k-1)$ fat blocks in any column of blocks $(S_{ij} : i \in [m'])$.
\begin{proof}
This is Lemma 5 of \cite{MT}. If there are more than ${t \choose {2k}}(k-1)$ fat blocks in a given column $S_{ij}$, then there are at least $k$ fat blocks which contain $1$'s in the same $2k$ columns of $A'$. Hence $A'$ contains a complete $k \times 2k$ matrix (i.e., a matrix in which all entries are $1$), so $A$ contains every $k \times 2k$ $(0,1)$--matrix, another contradiction.
\end{proof}

We wish to bound the number of $1$'s in $B$. $B$ may contain repeated rows, but since every row in $A'$ has at least $(2k-1)t + 1$ of its entries $1$, every row of $B$ must have at least $2k$ of its entries $1$. Thus if any row occurs in $B$ more than $k-1$ times, then $B$ contains a complete $k \times 2k$ matrix, contradicting Claim 1. If we let $B'$ be the matrix obtained from $B$ by deleting repeated rows, then $B'$ contains no matrix equivalent to $(I,M)$ and all rows of $B'$ are different, so at most $f(\lceil n/t \rceil,k)$ of the entries of $B'$ are $1$. Since each row in $B$ was repeated at most $(k-1)$ times, it follows that at most $(k-1)f(\lceil n/t \rceil,k)$ of the entries of $B$ are $1$. We have thus established the following recurrence: $$f(n,k) \le (2k-1)kf\left(\left\lceil \frac{n}{t} \right\rceil, k\right) + {t \choose {2k}} (k-1)n + g_D\left(2{D \choose 2}C_1(k)\right)n.$$

\noindent Let $t = 2k^2$. Assuming (by induction) that $$f \left( \left\lceil \frac{n}{2k^2} \right\rceil, k \right) < c_k \left\lceil \frac{n}{2k^2} \right\rceil,$$ and using the inequality $g_D(x) < 2D{x \choose {D-1}}$ noted earlier, we obtain
\begin{eqnarray}\label{eqn1}
f(n,k) & < & (2k^2-k)c_k \left\lceil \frac{n}{2k^2} \right\rceil + \left({{2k^2} \choose {2k}}k + 8k^3{{2{{4k^3 - 2k^2} \choose 2}C_1(k)} \choose {4k^3 - 2k^2 - 1}}\right)n \notag \\
& < & c_k n - \frac{c_k n}{2k} + (2k^2 - k)c_k + 8k^3{ {16k^6 C_1(k)} \choose {4k^3}}n.
\end{eqnarray}
So if
$$c_k \ge 32k^4{{16k^6C_1(k)} \choose {4k^3}},$$
then
\begin{equation}\label{eqn2}
\frac{c_k n}{2k} \: \ge \: \frac{c_k n}{4k} +  8k^3{{16k^6C_1(k)} \choose {4k^3}} \: \ge \: 2k^2c_k + 8k^3{{16k^6C_1(k)} \choose {4k^3}},
\end{equation} since $n \ge 8k^3$. Set $c_k = 32k^4{{16k^6C_1(k)} \choose {4k^3}}$, and note that $c_k > 2^{8k^3}$. Now inequalities (\ref{eqn1}) and (\ref{eqn2}) imply that $f(n,k) < c_kn$, so the induction step is complete.
\end{proof}

\begin{rmk}
Marcus and Tardos proved that $C(k) < 2k^4{{k^2} \choose k}$, so
the explicit bound we obtain is roughly $c_k = O(k^{9k^4})$.
\end{rmk}

\section{Ordered Hypergraphs}\label{hypersec}

We now deduce Theorem~\ref{hypergraphs} from Theorem~\ref{genMT}. This implication may be read out of a proof of Klazar (Theorem 2.5 of \cite{Klaz3}), but for the sake of completeness we shall prove it (and in fact our proof is slightly different from that in \cite{Klaz3}).

Given an ordered hypergraph $\HH$ on $[n]$, say that $\HH$ \emph{contains} a $k$-permutation $\pi$ if $\HH$ contains the ordered hypergraph $H(\pi)$. In other words, there exist $2k$ vertices $v_1, \ldots, v_{2k} \in [n]$, and $k$ distinct edges $E_1, \ldots, E_k \in E(\HH)$ such that, letting $e_i = \{v_i,v_{\pi(i)+k}\}$ denote the edges of $H(\pi)$, we have $e_i \subset E_i$ for each $i \in [k]$. Otherwise say that $\HH$ \emph{avoids} $\pi$. For each permutation $\pi$, and each $n \in \N$, let $T_n(\pi)$ denote the family of ordered hypergraphs on $[n]$ avoiding $\pi$.

\begin{lemma}\label{mattohg}
Let $k \in \N$ and $\pi \in \Pi_k$. If $\HH \in T_n(\pi)$, then $$\|\HH\| \: = \: \sum_{E \in E(\HH)} |E| \: < \: c_kn.$$
\end{lemma}

\begin{proof}
The lemma is a simple corollary of Theorem~\ref{genMT}. To spell it out, let $k \in \N$, $\pi \in \Pi_k$ and $\HH \in T_n(\pi)$, let $m = |E(\HH)|$, and define $A$ to be an $m \times n$ $(0,1)$--matrix whose rows are the indicator functions of the edges of $\HH$. (Note that $A$ is unique up to permutations of its rows.)

Now, $A$ has exactly $\|\HH\|$ of its entries 1, and each of its rows are different, so by Theorem~\ref{genMT}, if $\|\HH\| \ge c_kn$ then $A$ contains a sub-matrix $B \sim (I,M)$, where $M$ is the permutation matrix of $\pi$. Now, let $v_1, \ldots, v_{2k}$ be the vertices in $[n]$ corresponding to the columns of $B$, and $E_1, \ldots, E_k$ be the edges of $\HH$ corresponding to the rows of $B$, ordered so that $v_i \in E_i$ for each $i \in [k]$. Then for each $i \in [k]$ we have $v_i,v_{\pi(i)+k} \in E_i$, so $\HH$ contains $\pi$, a contradiction. Hence $\|\HH\| < c_kn$.
\end{proof}

\begin{proof}[Proof of Theorem~\ref{hypergraphs}]
Let $n,k \in \N$, and $\pi \in \Pi_k$. We claim that
\begin{equation}\label{eqn}
|T_{2n}(\pi)| \: \le \: |T_n(\pi)| \: 3^{2c^2n},
\end{equation}
where $c = c_k$ is the constant obtained in Theorem~\ref{genMT}.

To prove inequality~(\ref{eqn}), we map each $\HH \in T_{2n}(\pi)$ to the ordered hypergraph $\K$ on $[n]$ with edge set $$\{E \subset [n] : \exists E' \in E(\HH)\textup{ with }i \in E \Leftrightarrow \{2i-1,2i\} \cap E' \neq \emptyset\}.$$ In other words, $\K$ is formed by identifying vertices $2i-1$ and $2i$ for every $i \in [n]$.

We claim that $\K \in T_n(\pi)$, i.e., that $\K$ avoids $\pi$. Indeed, suppose for a contradiction that there exist edges $E_1, \ldots, E_k \in E(\K)$ and vertices $v_1, \ldots, v_{2k} \in [n]$ such that for each $i \in [k]$ we have $e_i = \{v_i,v_{\pi(i)+k}\} \subset E_i$. For each edge $E_i$ choose an edge $F_i$ of $\HH$ such that $\{2j-1,2j\} \cap F_i \neq \emptyset$ if and only if $j \in E_i$ (such an $F_i$ exists by the definition of $\K$). The edges $F_i$ are distinct (since the edges $E_i$ are), and $e_i \subset F_i$ for each $i \in [k]$, so $\HH$ contains $\pi$, which is the desired contradiction.

Now, how many ordered hypergraphs $\HH$ map to the same ordered hypergraph $\K$? An edge $E$ of $\K$ is the image of $3^{|E|}$ different possible edges of $\HH$, since each vertex $v$ of $E$ may have come from $\{2v-1\}$, $\{2v\}$ or $\{2v-1,2v\}$. However, suppose at least $2c = 2c_k$ of these did in fact occur in $\HH$ for a given edge $E$. Each such edge has size at least $|E|$, so the ordered hypergraph $\HH'$ induced by $\HH$ on vertex set $\{2i-1, 2i \in [2n] : i \in E\}$ has $\|\HH'\| \ge 2c|E|$. But now $\HH' \notin T_{2|E|}(\pi)$ by Lemma~\ref{mattohg}, so $\HH$ contains $\pi$, a contradiction.

Thus for each edge $E$ of $\K$, at most $2c - 1$ of the edges which map to it actually occur in $\HH$, so we have at most $$\sum_{i = 0}^{2c - 1} {{3^{|E|}} \choose i} < 3^{2c|E|}$$ choices for these edges.

Now, again by Lemma~\ref{mattohg}, since $\K \in T_n(\pi)$ we have $\|\K\| \le cn$. Thus the maximum possible number of ordered hypergraphs $\HH$ which map to a given $\K$ is
$$\prod_{E \in E(\K)} 3^{2c|E|} \: = \: 3^{2c \sum_E |E|} \: \le \: 3^{2c^2n}.$$ This proves inequality~(\ref{eqn}).

Now, let $\P$ be a strongly monotone property of ordered hypergraphs, let $k \in \N$, $\pi \in \Pi_k$, and suppose that $H(\pi) \notin \P$. Then $\P_n \subset T_n(\pi)$, and by inequality~(\ref{eqn}) and induction on $n$,
$$T_n(\pi) \: \le \: 3^{3c^2n},$$ for every $n \in \N$, where again $c = c_k$, since
$$|T_{2n-1}(\pi)| \: \le \: |T_{2n}(\pi)| \: \le \: |T_{n}| \: 3^{2c^2n} \: \le \: 3^{5c^2n} \: \le \: 3^{3c^2(2n - 1)}$$ when $n \ge 3$. So $|\P_n| \le 3^{3c^2n}$ for every $n \in \N$, which proves Theorem~\ref{hypergraphs}.
\end{proof}

\section{Partitions and ordered graphs}\label{rest}

We shall now deduce Theorems~\ref{parts}, \ref{mono} and \ref{noKt} from Theorem~\ref{hypergraphs}. We begin with Theorem~\ref{parts}. The implication is very simple, but in any case we shall write out all the details.

\begin{proof}[Proof of Theorem~\ref{parts}]
Let $\P$ be a hereditary property of partitions. For each partition $P \in \P$, let $\HH(P)$ be the ordered hypergraph whose edges are the parts of $P$ of size at least two. To be precise, if $P$ is the partition $\{A_1, \ldots, A_t\}$ of $[n]$, then $\HH(P)$ has vertex set $[n]$ and edge set $\{A_i : i \in [t], |A_i| \ge 2\}$. Observe that for any permutation $\pi$, $\HH(P)$ contains $\pi$ if and only if $P$ contains the partition $H(\pi)$ as an induced subpartition.

Suppose that for some $\pi \in \Pi$, $\P$ does not contain $H(\pi)$. By the observation above, $\HH(P)$ avoids $\pi$ for every $P \in \P$. Now, let $T(\pi)$ denote the strongly monotone property of ordered hypergraphs consisting of all ordered hypergraphs avoiding $\pi$ (so $T(\pi) = \bigcup_n T_n(\pi)$, with $T_n(\pi)$ as in the previous section). Then $\HH(P) \in T(\pi)$ for every $P \in \P$.

Now apply Theorem~\ref{hypergraphs} to $T(\pi)$. Since $H(\pi) \notin T(\pi)$, there exists a constant $c$ such that $|T_n(\pi)| \le c^n$ for every $n \in \N$. But now we are done, since
$$|\P_n|\: = \:|\{\HH(P) : P \in \P_n\}| \: \le \: |T_n(\pi)| \: \le \: c^n$$
for every $n \in \N$, since $\HH(P) \in T(\pi)$ for every $P \in \P$. This proves Theorem~\ref{parts}.
\end{proof}

\begin{rmk}
If $\pi \in \Pi$, let $\P(\pi)$ be the largest hereditary property of partitions such that $\P$ avoids the partition $H(\pi)$. Let $c'_k$ be the smallest constant such that $|\P(\pi)_n| < (c_k')^n$ for every $n \in \N$ and $\pi \in \Pi_k$. We have shown that $c'_k = O(3^{k^{19k^4}})$.
\end{rmk}

We next deduce Theorem~\ref{mono} from Theorem~\ref{hypergraphs}. Since a monotone property of ordered graphs is a strongly monotone property of ordered hypergraphs, the implication is trivial.

\begin{proof}[First proof of Theorem~\ref{mono}]
Let $\P$ be a monotone property of ordered graphs, and for each $G \in \P$, let $G'$ be the ordered hypergraph with the same vertex and edge set as $G$. Define $\P' = \{G' : G \in \P\}$. Now $\P'$ is a strongly monotone property of ordered hypergraphs, since each edge of $G' \in \P'$ has size 2, so the only ordered hypergraphs contained in $G'$ are its subgraphs. The result now follows by applying Theorem~\ref{hypergraphs} to $\P'$.
\end{proof}

In fact one can also prove Theorem~\ref{mono} without using Theorem~\ref{hypergraphs}, but using the Marcus-Tardos and Klazar-Marcus-Tardos Theorems instead. Alternative proofs can often give new insight into the difficulties and the true nature of a problem, and for this reason we give a sketch of this second proof.

\begin{proof}[Sketch of the second proof of Theorem~\ref{mono}]
Let $\P$ be a monotone property of ordered graphs, let $\pi \in \Pi$, and suppose that $\P$ does not contain the ordered graph $H(\pi)$.

Suppose first that for some $n \in \N$ there exists an ordered graph $G \in \P_n$ with at least $C(k+1)n$ edges (where $C(k)$ is again the constant in Theorem B). In this case we can use Theorem B to find $H(\pi)$ in $G$, just as in the proof of Lemma~\ref{MTdeg2}. So assume that for every $G \in \P$, $e(G) < C(k+1)|G|$.

Let $\S(n,m)$ denote the family of sequences $(a_1, \ldots, a_n)$ such that $a_i \in \N \cup \{0\}$ for each $i$ and $\sum_i a_i = m$, and let $\S = \bigcup_{n,m} S(n,m)$. We define a map $\varphi : \P \to \Pi \times \S \times \S$ as follows.

Let $G \in \P$ have $n$ vertices and $m$ edges. We put the following two linear orders, $<_\ell$ and $<_r$, on the edges of $G$. If $e = \{e_1,e_2\}$ and $f = \{f_1,f_2\}$ with $e_1 < e_2$ and $f_1 < f_2$, then $e <_\ell f$ if $e_1 < f_1$, or $e_1 = f_1$ and $e_2 < f_2$, while $e <_r f$ if $e_2 < f_2$, or $e_2 = f_2$ and $e_1 < f_1$. Let $\varphi_p(G)$ be the $m$-permutation which takes the order of the edges under $<_r$ to the order under $<_\ell$. Let $\varphi_\ell (G)$ be the left-endpoint degree sequence of $G$, i.e., the sequence $(a_1, \ldots, a_n)$ where $a_i$ is the number of edges of $G$ whose left-endpoint is vertex $i$, and let $\varphi_r (G)$ be the right-endpoint degree sequence of $G$. Let $\varphi(G) = (\varphi_p(G), \varphi_\ell (G), \varphi_r(G)) \in \Pi_m \times S(n,m) \times S(n,m)$.

Let $\Q = \{\varphi_p(G) : G \in \P\}$. $\Q$ is a hereditary property of permutations, so by the Klazar-Marcus-Tardos Theorem, either $Q = \Pi$, or there exists a constant $c$ such that $|\Q_n| \le c^n$ for every $n \in \N$.

Suppose first that $|\Q_n| \le c^n$ for every $n \in \N$. We claim that for any $\varphi \in \Pi \times \S \times \S$, there is only at most one $G$ such that $\varphi(G) = \varphi$. We omit the proof, which is by induction on $m$. For the induction step, remove the first edge of $G$ in the order $<_r$.

So $|\P_n|$ is just $|\operatorname{Im}_{\varphi}(\P_n)|$, which can easily be approximated since each $G \in \P_n$ has at most $C(k+1)n$ edges. Thus
\begin{eqnarray*}
|\operatorname{Im}_{\varphi_p}(\P_n)| & \le & |\bigcup_{m=0}^{C(k+1)n} \Q_m|\\
& \le & \sum_{k=0}^{C(k+1)n} c^m \: < \; 2c^{C(k+1)n},
\end{eqnarray*}
assuming (as we may) that $c > 2$, and
$$|\operatorname{Im}_{\varphi_\ell}(\P_n)| \le {{m+n-1} \choose {n-1}} < {{(C(k+1)+1)n} \choose {n-1}} < 2^{(C(k+1)+1)n},$$ and similarly for $|\operatorname{Im}_{\varphi_r}(\P_n)|$. Hence
\begin{eqnarray*}
|\P_n| = |\operatorname{Im}_{\varphi}(\P_n)| & \le & |\operatorname{Im}_{\varphi_p}(\P_n)| \cdot |\operatorname{Im}_{\varphi_\ell}(\P_n)| \cdot |\operatorname{Im}_{\varphi_r}(\P_n)|\\[+1ex] & < & 2^{2(C(k+1)+1)n+1}c^{C(k+1)n},
\end{eqnarray*}
so we are done in this case.

Now suppose that $\Q = \Pi$. We want to show that $\P$ contains the ordered graph $H(\pi)$, and thus obtain a contradiction. To do this, define the $2k$-permutation $\sigma$ by $\sigma(2i-1) = 2\pi(i)$ and $\sigma(2i) = 2\pi(i)-1$ for $1 \le i \le k$ (so for example if $\pi = 213$ then $\sigma = 432165$). By assumption, there exists an ordered graph $G \in \P$ such that $\varphi_p(G) = \sigma$. Let the edges of $G$ be $e_1, \ldots, e_{2k}$ in the order $<_\ell$. Note that for each $i \in [k]$, the edges $e_{2i-1}$ and $e_{2i}$ do not share an endpoint, since $e_{2i-1} <_\ell e_{2i}$ and $e_{2i-1} >_r e_{2i}$.

Consider the edges $\{e_{2i-1} : i \in [k]\}$. Suppose two of them share a left (right) endpoint $v$. Then all edges between them in the order $<_\ell$ ($<_r$) also share that endpoint. This contradicts the previous observation that $e_{2i-1}$ and $e_{2i}$ do not share an endpoint, so these edges are in fact independent. Thus there exists $G \in \P$ with $\Delta(G) = 1$ and $\varphi_p(G) = \pi$.

Now apply the same technique to the $(k+1)$-permutation $\pi'$, where $\pi'(i) = \pi(i) + 1$ if $i \in [k]$, and $\pi'(k+1) = 1$. We obtain $G \in \P$ with $\Delta(G) = 1$ and $\varphi_p(G) = \pi'$. Again (as in the proof of Lemma~\ref{deg1}), notice that in $G$ all left-endpoints occur to the left of \emph{all} right-endpoints, since $\pi'(k+1) = 1$. Therefore, letting $G'$ be subgraph of $G$ induced by the first $k$ left-endpoints and last $k$ right-endpoints, we have $G' = H(\pi)$.

So $H(\pi) \in \P$, and this gives us the desired contradiction.
\end{proof}

Finally, we deduce Theorem~\ref{noKt} from Theorem~\ref{mono}.

\begin{proof}
Let $t \in \N$, and let $\P$ be a hereditary property of ordered graphs such that $K_t \notin \P$ and $K_{t,t} \notin \P$. Let $\G = \G(\P)$ be the smallest monotone property containing $\P$. Note that (trivially) $|\G_n| \ge |\P_n|$.

By Theorem~\ref{mono}, either $\G$ contains $H(\pi)$ for every $\pi \in \Pi$, or there exists a constant $c$ such that $|\G_n| < c^n$ for every $n \in \N$. Suppose the latter. Then $|\P_n| \le |\G_n| < c^n$ for every $n \in \N$, in which case we are done. So assume that $\G$ contains $H(\pi)$ for every $\pi \in \Pi$.

We wish to find, for each permutation $\pi$, a large permutation $\sigma$ such that an ordered graph on $[n]$ containing no induced copy of $K_t$ or $K_{t,t}$, and containing $H(\sigma)$ as a subgraph, contains $H(\pi)$ as an induced subgraph. We shall use Ramsey's Theorem and the pigeonhole principle, so recall that for $\ell \in \N$, $R(\ell)$ denotes the smallest integer $n$ such that any graph on $n$ vertices contains either a clique on $\ell$ vertices, or an independent set of order $\ell$. We shall also write $S(\ell)$ for the smallest integer $n$ such that any bipartite graph with one part of order at least $2\ell - 1$ and the other of order at least $n$, contains either the complete bipartite graph $K_{\ell,\ell}$, or the empty bipartite graph $E_{\ell,\ell}$. It is easy to show that $R(\ell) < 4^\ell$ and $S(\ell) < \ell4^\ell$. For each $j \in \N$, let $R^{(j+1)}(\ell) = R(R^{(j)}(\ell))$, where $R^{(1)}(\ell) = R(\ell)$, and let $S^{(j)}(\ell)$ be defined similarly.

Let $k \in \N$ and $\pi \in \Pi_k$. We shall show that $H(\pi) \in \P$. Let $m' = S^{(K)}(t)$, where $K = {{2k} \choose 2} - k$, and $m = R^{(2)}(m')$. We define the $mk$-permutation $\sigma$ as follows: for each $i \in [k]$ and $j \in [0,m-1]$, let $\sigma(im - j) = \pi(i)m - j$. By assumption, $H(\sigma) \in \G$; let $G \in \P$ be an ordered graph on $[2mk]$ containing $H(\sigma)$. Such a $G$ must exist by the definition of $\G$. We know that for each $i \in [mk]$, the edge $\{i,\sigma(i) + mk\}$ is in $E(G)$. We want to show that for some subset $A = \{a(1), \ldots, a(2k)\} \subset [2mk]$ with $a(i) \in [(i-1)m + 1, im]$ for each $i \in [2k]$, and $a(\pi(i)+k) = a(i) + (k + \pi(i) - i)m$ for each $i \in [k]$, these are the only edges induced by $A$.

First we use Ramsey's Theorem to find `matching' independent subsets of $[(i-1)m + 1,im]$ for each $i \in [2k]$. (This step is necessary if we are to assume only that $K_t$ and $K_{t,t}$ are avoided; without it we would have to assume that every ordered graph containing $K_{t,t}$ is missing from $\P$.) We do this first for each $i \in [k]$. Since $m \ge R^{(2)}(m')$, by Ramsey's Theorem there exists either a clique or an independent set $A^{(1)}_i$ of order $R(m')$ in $[(i-1)m + 1,im]$. Since, $R(m') \ge t$, $K_t \notin \P$ and $\P$ is hereditary, $A^{(1)}_i$ must be an independent set. For $i \in [k+1,2k]$, let $A^{(1)}_i$ be the set $A^{(1)}_{\pi^{-1}(i-k)} + (k + \pi(i) - i)m$ (if $A$ is a set and $b \in \N$ then $A + b = \{a + b : a \in A\}$). Now, again by Ramsey's Theorem, for each $i \in [k+1,2k]$ there exists a clique or independent set $A^{(2)}_i \subset A^{(1)}_i$ of order $m'$. Again this must be an independent set, since $m' \ge t$. For each $i \in [k]$, let $A^{(2)}_i$ be the set $A^{(2)}_{\pi(i)+k} - (k + \pi(i) - i)m$.

Thus we have found independent sets $A^{(2)}_i \subset [(i-1)m+1,im]$ of order $m'$ for each $i \in [2k]$ such that there is a matching between $A^{(2)}_i$ and $A^{(2)}_{\pi(i)+k}$ for each $i \in [k]$. We next apply the pigeonhole principle (aka bipartite Ramsey Theorem) to each of the pairs $A^{(2)}_x$ and $A^{(2)}_y$ with $1 \le x < y \le 2k$, $y \neq \pi(x) + k$, to find the desired subset $A$.

To be precise, let $\{e_1, \ldots, e_K\}$ be the set of pairs $\{x,y\}$ such that $1 \le x < y \le 2k$ and $y \neq \pi(x) + k$ (recall that $K = {{2k} \choose 2} - k$), and for each $i \in [2k]$ let $B^{(0)}_i = A^{(2)}_i$. We shall define inductively, for each $i \in [2k]$, a sequence of sets $B^{(0)}_i \supset B^{(1)}_i \supset \ldots \supset B^{(K)}_i$. For each pair $e_\ell = \{x(\ell), y(\ell)\}$ in turn (i.e., for each $\ell \in [K]$), define the sets $B^{(\ell)}_i$ as follows.

Let $x' = x'(\ell)$ and $y' = y'(\ell)$ be the elements matched to $x = x(\ell)$ and $y = y(\ell)$ respectively by $\pi$, so $x' = \pi(x) + k$ if $x \le k$ and $x' = \pi^{-1}(x-k)$ if $x \ge k+1$, and similarly for $y'$. If $i \notin \{x, y, x', y'\}$ then set $B^{(\ell)}_i = B^{(\ell - 1)}_i$. Let $B^{(\ell)}_{x}$ and $B^{(\ell)}_{y}$ be the parts of the largest empty bipartite graph induced by $G$ with $B^{(\ell)}_{x} \subset B^{(\ell-1)}_{x}$, $B^{(\ell)}_{y} \subset B^{(\ell-1)}_{y}$ and $|B^{(\ell)}_{x}| = |B^{(\ell)}_{y}|$. Let $B^{(\ell)}_{x'}$ be the set $B^{(\ell)}_{x} + (k + \pi(i) - i)m$ if $x \le k$ and the set $B^{(\ell)}_{x} - (k + \pi(i) - i)m$ if $x \ge k + 1$. Let $B^{(\ell)}_{y'}$ be defined from $B^{(\ell)}_{y}$ similarly. Note that $B^{(\ell)}_{x'} \subset B^{(\ell-1)}_{x'}$, and $B^{(\ell)}_{y'} \subset B^{(\ell-1)}_{y'}$.

We claim that $|B^{(j)}_i| \ge S^{(K - j)}(t)$ for each $i \in [2k]$ and $j \in [K]$, and prove it by induction on $j$. For $j = 0$ the statement is that $|B^{(0)}_i| \ge S^{(K)}(t) = m'$ for each $i \in [k]$, so the base case holds (since each set $A^{(2)}_i$ has order $m'$). Assume the result holds for $j-1$. If $i \notin \{x, y, x', y'\}$, then $|B^{(j)}_i|  = |B^{(j-1)}_i| \ge S^{(K - j + 1)}(t) > S^{(K - j)}(t)$, so we are done in this case. Now suppose $i \in \{x, y, x', y'\}$. By the induction hypothesis, $|B^{(j-1)}_\ell| \ge S^{(K - j + 1)}(t)$ for each $\ell \in \{x, y\}$, so by the definition of $S$, there exists either a complete bipartite or an empty bipartite graph in $G[B^{(j-1)}_{x},B^{(j-1)}_{y}]$ with each part having at least $S^{(K - j)}(t)$ vertices. Since $S^{(K - j)}(t) \ge t$, it cannot be complete (since the sets $B^{(j-1)}_{x}$ and $B^{(j-1)}_{y}$ are independent), so $B^{(j)}_\ell \ge S^{(K - j)}(t)$, for $\ell \in \{x,y\}$. Since $|B^{(j)}_{x'}| = |B^{(j)}_x|$ and $|B^{(j)}_{y'}| = |B^{(j)}_y|$, the induction step is complete.

It follows from the claim that $|B^{(K)}_i| \ge t \ge 1$ for each $i \in [2k]$. Observe also that for each $i \in [k]$, $B^{(K)}_i + (k + \pi(i) - i)m = B^{(K)}_{\pi(i)+k}$, and that for $\{i,j\} \in \{e_1, \ldots, e_K\}$, there are no edges in $G$ between $B^{(K)}_i$ and $B^{(K)}_j$. For each $i \in [k]$, choose a vertex $a(i) \in B^{(K)}_i$, and let $a(\pi(i)+k) = a(i) + (k + \pi(i) - i)m \in B^{(K)}_{\pi(i)+k}$. Let $A = \{a(1), \ldots, a(2k)\}$. The set $A$ induces the ordered graph $H(\pi)$ in $G$, so we are done.
\end{proof}

\begin{rmk}
For each permutation $\pi$, let $\tilde{c}(\pi, t)$ denote the smallest constant such that $|\P_n| < \tilde{c}(\pi, t)^n$ for every $n \in \N$, for every hereditary property $\P$, satisfying the conditions of the theorem, which avoids $H(\pi)$. The bounds given by our proof on the constant $\tilde{c}(\pi,t)$ are rather large. They could be improved somewhat by choosing the order in which the pairs $e_i$ are dealt with, and thus obtaining a much stronger inequality than the one we obtained ($|B^{(j)}_i| \ge S^{(K - j)}(t)$), but for simplicity of presentation (and because the actual bounds are not our main interest), we leave this as an exercise for the interested reader. Notice also that although we assumed $K_t \notin \P$, we only needed $K_{m'} \notin \P$.
\end{rmk}

We finish by noting an immediate consequence of Theorem~\ref{noKt}.

\begin{cor}
Let $\P$ be a hereditary property of ordered graphs. If there exists a function $f: \N \to \N$ such that $$e(G) \le f(n) = o(n^2)$$ for every $G \in \P_n$ and every $n \in \N$, then Conjecture~\ref{orderconj} holds for $\P$.
\end{cor}

\begin{proof}
Suppose there is such a function $f(n) = o(n^2)$, satisfying $e(G) \le f(n)$ for every $G \in \P_n$ and every $n \in \N$. Since $f(n) = o(n^2)$, there must exist $t \in \N$ such that $e(G) < n^2/4$ for every $G \in \P_n$ with $n \ge t$.

Assume $t \ge 2$. Now, $e(K_t) = {t \choose 2} \ge \frac{t^2}{4}$, and $e(K_{t,t}) = t^2$, so $K_t \notin \P_t$ and $K_{t,t} \notin \P_{2t}$. The result now follows by Theorem~\ref{noKt}.
\end{proof}

\section{Acknowledgements}

The authors would like to thank the anonymous referees for their careful reading of the manuscript and their many helpful comments, which included simplifying the original proof of Lemma~\ref{bound}.

\end{document}